\documentclass{gtart}

\usepackage[dvips]{graphicx}

\usepackage{pinlabel}

\usepackage{amssymb,amsfonts,amsmath,amsthm}

\newtheorem{theorem}{Theorem}

\theoremstyle{definition}

\begin{document}

\sloppy



\title{Fiber connected, indefinite Morse $2$-functions on connected
$n$-manifolds.
}





%
%

\authors{David T. Gay, Robion Kirby
\footnote{Supported in part by South African National Research Foundation Focus Area grant FA2007042500033 and United States National Science Foundation grant EMSW21-RTG} 
}
\address{Euclid Lab, 428 Kimball Rd, Iowa City, IA 52245\\
Department of Mathematics,University of Iowa, Iowa City, IA 52242}
\secondaddress{University of California, Berkeley, CA 94720} 
\email{d.gay@euclidlab.org}
\secondemail{kirby@math.berkeley.edu}



\begin{abstract} We discuss generic smooth maps from smooth manifolds to smooth surfaces, which we call ``Morse $2$--functions'', and homotopies between such maps. The two central issues are to keep the fibers connected, in which case the Morse $2$--function is ``fiber-connected'', and to avoid local extrema over $1$--dimensional submanifolds of the range, in which case the Morse $2$--function is ``indefinite''. This is foundational work for the long-range goal of defining smooth invariants from Morse $2$--functions using tools analogous to classical Morse homology and Cerf theory. \end{abstract}







\maketitlepage

Let $X^n$ be a smooth ($C^{\infty}$), compact, oriented, connected,
$n$-manifold. We consider generic, smooth maps to a $1-$ or
$2$-dimensional manifold and state theorems about their critical
points, whether they are indefinite or not, whether level sets or
fibers are connected, and to what extent these properties hold in
$1$-parameter families.

The classical case is that of real valued Morse functions, $f:S^n \to
R^1$, where the differential $df$ either has rank one, or if the rank
is zero then the Jacobian matrix $(\partial^2f/\partial x_i\partial
x_j)$ is non-singular.  Since $X$ is connected, it is well known that
$f$ need have only one minimum (equal definite critical point of index
0) and one maximum (definite of index $n$) and that level sets are
connected (for $n \geq 3$ because $f$ can be chosen to be
self-indexing.  (Self indexing means that if two critical points, $x$
and $Y$ have index $i$ and $j$ respectively, with $i < j$, then $f(x)
< f(y)$;  it follows that passing a critical point of index $n-1$ will
not disconnect the level set because  $X$ is connected and $n \geq 3$.)

The relative case $f:(X,\partial X) \to (B^1, S^0)$ is similar, but
note that there need be no definite critical points.

Of more interest here are circle valued Morse functions on closed
$n$-manifolds, $f:S^n \to S^1$.  We assume that the induced map on
first cohomology, $f^* : Z = H^1(S^1;Z) \to H^1(X;Z)$ sends the
generator $1 \in Z$ to a primitive class in $H^1(X;Z)$.  (This avoids,
for example, the map $Y^{n-1} \times S^1 \to S^1$ sending $(y,\theta)
\to 2\theta$ with disconnected fibers.)  Then, $f$ is homotopic
to a Morse function without definite critical points and with
connected fibers. (For maps from manifolds to manifolds, point-inverses are usually called ``fibers'', except in the classical case of a map to $R$, which is often thought of as a height function, and point-inverses are called level sets.)

Now suppose we are given two homotopic Morse functions $f_0, f_1:(X^n,
\partial X) \to (B^1, S^0)$ which are indefinite with connected level
sets. Then Cerf theory \cite{Cerf} tells us that the homotopy, $f_t$, can be
chosen to consist of Morse functions except for finitely many values
of $t$ at which either births or deaths of canceling critical points
occur.  The local model for such is 
$$f_t (u,y) = u^3 -tu -y_1^2 -\cdots -y_k^2 +y_{k+1}^2 + \cdots +
y_n^2$$ for $(u,y) \in R \times R^{n-1}$, where there are no critical
points for $t<0$, a birth at $t=0$, and critical points of index $k$
and $k+1$ for $t>0$. It was shown in \cite{Kirby} that $f_t$ may also
be chosen so that there are no births of definite critical points and
level sets always remain connected.

\begin{figure} 
\labellist
\tiny\hair 2pt
\pinlabel $1$ [r] at 6 25
\pinlabel $2$ [r] at 6 46
\pinlabel $2$ [r] at 6 67
\pinlabel $3$ [r] at 6 88
\pinlabel $1$ [l] at 210 25
\pinlabel $2$ [l] at 210 61
\pinlabel $1$ [tr] at 37 34
\pinlabel $2$ [br] at 32 38
\pinlabel $2$ [tr] at 50 72
\pinlabel $3$ [br] at 45 80
\pinlabel $1$ [br] at 94 38

\endlabellist

\begin{center}

\includegraphics{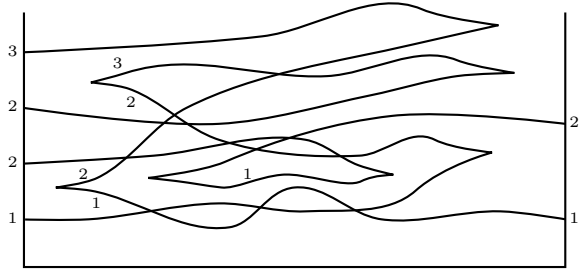}
\caption{A typical Cerf graphic.} \label{F:CerfGraphic}
 
\end{center}
\end{figure}

The term {\it Cerf graphic} is used for the arcs in $I \times B^1$
which are points $(t,v) \in I \times B^1$ for which $v$ is a critical
value for $f_t$.  An example is given in Figure~\ref{F:CerfGraphic} for $n=4$, where
there are only critical points of index $1, 2$, and $3$.  Note that
births can occur earlier, deaths later, and $1$'s can be pushed below
$2$'s, etc.  Not illustrated in Cerf graphics are the times $t$ at
which, for example, the descending $1$-manifold of a critical point of index $1$
hits a lower critical point of index $1$;  these correspond to a
$1$-handle sliding over another $1$-handle.

Finally, consider homotopic maps $f_0, f_1 :S^n \to S^1$ with
connected fibers and no definite critical points.

\begin{theorem}[\cite{GayKirbyMorse2Fs}]  The homotopy $f_t$ between $f_0$ and $f_1$ is
itself homotopic, fixing $f_0$ and $f_1$, to a homotopy still called
$f_t$ which has no births or deaths of definite critical points and
during which the fibers remain connected.
\end{theorem}

It is easy to prove that definite births and deaths need not occur
(the same argument is used as for range equal to $R^1$), but it takes
some work to keep fibers connected because it is not trivial to make
the $f_t, t\in (0,1)$ self indexing.  We thank Katrin Wehrheim and Chris Woodward for
bringing this issue to our attention.

The main point of this paper is to address the cases $p:X^n \to
\Sigma^2$ where $p^{-1}(\partial \Sigma) = \partial X$, and $\Sigma$
is the $2$-ball $B^2$ (sometimes conveniently the square $I \times
I$), or the $2$-sphere $S^2$, or the annulus $I \times S^1$, or a
general oriented, connected, compact surface.

According to singularity theory~\cite{Wassermann, LekiliWrinkled} the rank of the differential of a generic
smooth  map $p:X \to \Sigma$ will be two except along a smooth,
embedded $1$-manifold $\hat{L}$ in $X$ which maps by $p$ to an
immersed $1$-manifold $L$ in $\Sigma$ with cusps (occurring when
$\hat{L}$ has vertical tangents).

\begin{figure} 
\labellist
\tiny\hair 2pt
\pinlabel $1$ [l] at 84 41
\pinlabel $2$ [tr] at 56 79

\endlabellist

\begin{center}

\includegraphics{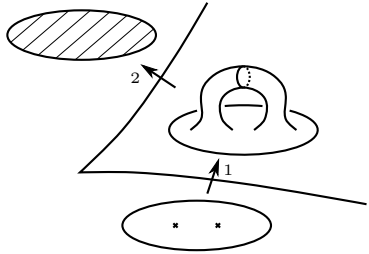}
\caption{Fibers near a cusp and two folds.} \label{F:CuspAndFold}
\end{center}
\end{figure}

Locally the map $p$ is modeled by pieces of the Cerf graphic for $f$
as described above.  A non-cusp point of $\hat{L}$ has a neighborhood
with local coordinates $(t,y) \in R \times R^{n-1}$ such that $p(t,y)
= (t,\sum \pm y^2_i)$ so that $p$ has an arc of critical points $(t,0)
\subset \hat{L}$.  If the critical points have index $k$, then locally
$X$ is $R$ cross the standard $(n-1)$-bordism in which a $k$-handle
(or $(n-1-k)$-handle from the opposite direction) is attached to the
fiber $F_1$ on one side of $L$ to get the fiber $F_2$ on the other
side, see Figure~\ref{F:CuspAndFold}. $L$ is called a {\it fold} or more precisely a
$k$-fold or $(n-1-k)$-fold.

At a cusp, there are local coordinates $(t,y) \in R \times R^{n-1}$
for which $p(t,y) = (t,y_2^3 -ty_2 + \sum_{i=3}^n  \pm y_i)$.

 We call such generic maps {\it Morse $2$-functions}, to emphasize the connection with classical Morse theory and to emphasize that the target space is $2$--dimensional.   As with the
case of Morse functions $f:X \to B^1$ or $S^1$ where $X$ is connected,
we are interested in when $p:X \to \Sigma$ can be homotoped to be
indefinite (this is analogous to no maxima or minima) and to have
connected fibers (analogous to connected level sets).  The answer is
the following existence theorem.

\begin{theorem}[Existence] \label{T:Existence}
Let $p:(X, \partial X) \to (\Sigma, \partial
\Sigma)$ be a map which restricted to $\partial X$ is an indefinite,
surjective Morse function.  If $n>2$ and $p_*(\pi_1 (X))$ has finite
index in $\pi_1 (\Sigma)$, the $p$ is homotopic rel boundary to an
indefinite Morse $2$-function, still called $p:X \to \Sigma$.  If $n
>3$, $p$ restricted to $\partial \Sigma$ is fiber-connected, and if
$p_*(\pi_1 (X))= \pi_1 (\Sigma)$, then we can choose $p$ to be fiber
connected.
\end{theorem}

The indefinite part of this theorem was proved by Osamu Saeki
\cite{Saeki} for $\Sigma = S^2$ or $RP^2$, and the necessary
conditions on $\pi_1$ were mentioned.  This result for $4$--manifolds mapping to $S^2$, with fiber-connectedness, follows from earlier work of the authors~\cite{GayKirbyBLFs} and the work of Lekili~\cite{LekiliWrinkled}, Baykur~\cite{BaykurExist} and Akbulut and Karakurt~\cite{AkbulutKarakurt}. 

What data, besides the folds in the image $p(X)$ in $S^2$, determines
$X$ up to diffeomorphism?  In the case of a Morse function $f:X \to
R$, one needs to know the descending manifold for each critical point
of $f$, or equivalently, the attaching framed $(k-1)$-sphere in the
level set just below a critical point of index $k$. This is the local data; globally, all the data of the attaching maps can be drawn in the boundary of the $0$--handle.

A similar statement is true
for $p:X^4 \to S^2$. We begin with the local picture: Because $n=4$, there are only $1$- and
$2$-folds and these are not actually different since a $1$-fold from
one direction is a $2$-fold from the other.  The fiber $F$ is
$2$-dimensional and crossing a $2$-fold means adding a $2$-handle to a
non-separating circle in $F$. In the opposite direction, crossing a $1$--fold means adding a $1$--handle to a pair of points in $F$. These circles and pairs of points are part of the local data.  

\begin{figure} 
\labellist
\tiny\hair 2pt
\pinlabel $2$ [br] at 73 198
\pinlabel $2$ [tr] at 77 158
\pinlabel $a$ [b] at 109 203
\pinlabel $b$ [b] at 110 175
\small
\pinlabel (a) at 51 153
\pinlabel (b) at 45 34

\endlabellist

\begin{center}

\includegraphics{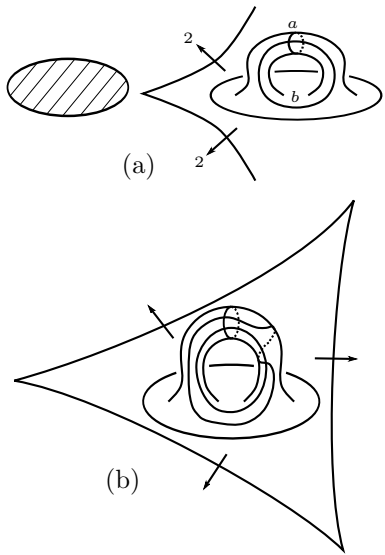}
\caption{Cusps and Lefschetz singularities.} \label{F:CuspLefschetz}
\end{center}

\end{figure}

At a crossing, we have four surfaces, one of genus $g+2$, two of genus $g+1$ and one of genus $g$. The attaching circles in the genus $g+2$ surface must be disjoint and must be identified with the attaching circles in the genus $g+1$ surfaces in the obvious way. At a cusp, the two attaching circles on the high genus side must meet in exactly one point, see Figure~\ref{F:CuspLefschetz}a. (In~\cite{LekiliWrinkled}, Lekili showed that the cusp may be smoothed if it is replaced by a Lefschetz singularity with vanishing cycle equal to a circle $C$ such that a Dehn twist about $C$ carries $a$ to $b$, as in Figure~\ref{F:CuspLefschetz}a.  Removing a Lefschetz singularity introduces the triangle of folds in Figure~\ref{F:CuspLefschetz}b.)

\begin{theorem}
 This local data determines the $4$--manifold if all fibers have genus greater than $1$ and all regions bounded by folds in $S^2$ are simply connected.
\end{theorem}

\textbf{Proof:} 
Let $F \subset S^2$ be the image of the fold set, and let $F'$ be a
parallel copy of $F$ pushed off in the direction of decreasing genus.
Note that, for each polygonal region $R$ in $S^2 - F'$, the local data
completely determines $p^{-1}(R)$. The vertices of $\partial R$ occur near crossings of $F$; at each crossing there are four quadrants, with fibers of genus $g+2$, $g+1$, $g+1$ and $g$, and there will be a vertex in the genus $g$ quadrant. Consider two adjacent regions $R_1$
and $R_2$, meeting along an edge $e$. There will be one or two arcs $a \subset
F$ crossing $e$. For each such arc, we first glue $p^{-1}(R_1)$ and $p^{-1}(R_2)$
together along a small neighborhood of $a$; there may be many choices
for how to do this, but the local data should patch together
consistently. Having made these choices, what remains to be filled in
is the inverse image of a collection of disks with no singularities;
in other words, we need to glue in $F \times D^2$ for various fibers
$F$. The fact that the original local data came from a map $p : X^4
\to S^2$ means that there exist choices such that the boundary
monodromies are trivial and thus so that $F \times D^2$ can be glued
in. If we make different choices which still yield trivial boundary
monodromies, we may end up gluing $F \times D^2$ in differently, but
the fact that the genus of each fiber is greater than $1$ means that
$\pi_1(\mathit{Diffeo}_+(F)) = 0$, which means that different ways of
gluing in $F \times D^2$ give the same $4$--manifold.
$\Box$

%

%
%
%
%

If some region bounded by folds is not simply connected, then nontrivial monodromy may appear and must be specified. If the
genus of $F$ is one in some region, then a nontrivial element of $\pi_1(\mathit{Diffeo}_+(F))$ may appear and must also be
specified. An interesting example of this appears in the end of this
paper. 

When $F = S^2$, there is one non-trivial loop of diffeomorphisms,
namely rotate $S^2$ by $\theta$ for $\theta$ in the loop.  Any knotted
$2$-sphere in $S^4$ can be realized as the fiber $F$ of a
fiber-connected, indefinite Morse $2$-function  $S^4 \to
S^2$ (which is necessarily homotopic to the constant map).  Then this
$2$-sphere fiber, $F$, can be cut out and glued back in by the
non-trivial loop above.  This is called the {\it Gluck twist}
\cite{Gluck} and it is very interesting because the resulting $X$ is
homotopy equivalent to $S^4$, therefore homeomorphic to $S^4$
\cite{Freedman}, but not known to be diffeomorphic to $S^4$ in general.  Whether or not a smooth homotopy $4$--sphere is diffeomorphic to $S^4$ is
the last remaining case of the famous Poincar\'{e} Conjectures.

\begin{figure} 
\labellist
\tiny\hair 2pt
\pinlabel $k$ [b] at 50 180
\pinlabel $k$ [b] at 190 192
\pinlabel $k+1$ [b] at 161 204
\pinlabel $k$ [b] at 190 168
\pinlabel $k-1$ [t] at 162 156
\pinlabel $k$ [b] at 160 96
\pinlabel $k+1$ [b] at 160 120
\pinlabel $k$ [tl] at 28 32
\pinlabel $k+1$ [bl] at 28 41
\pinlabel $k$ [tr] at 76 32
\pinlabel $k+1$ [br] at 76 41
\pinlabel $k$ [t] at 160 30
\pinlabel $k+1$ [b] at 160 42

\small
\pinlabel (a) [r] at 16 180
\pinlabel (b) [r] at 16 108
\pinlabel (c) [r] at 16 38
\pinlabel ? [b] at 107 29

\large
\pinlabel $\emptyset$ at 57 108

\endlabellist

\begin{center}

\includegraphics{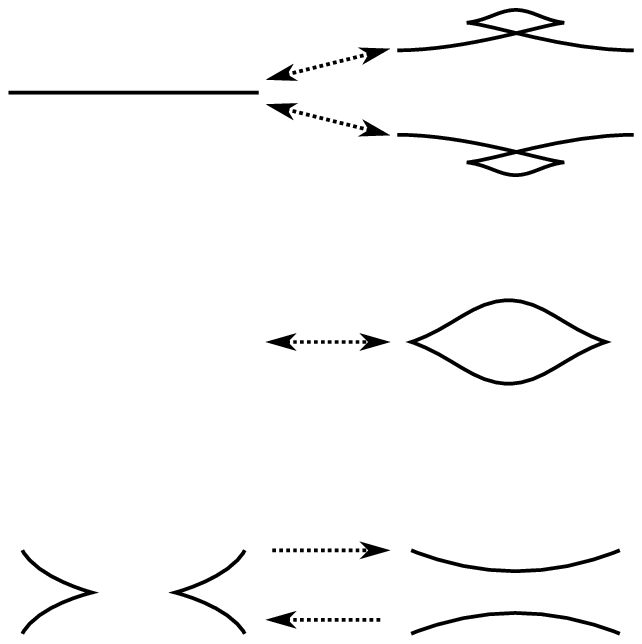}
\caption{Moves on a Cerf graphic.} \label{F:Moves}
\end{center}
\end{figure}

Before giving a proof of existence for the case $\Sigma= S^2$, we
address the question of uniqueness.  From singularity theory again,
there are essentially three local changes that can be made ( besides
regular homotopies, called ``isotopies'' in~\cite{LekiliWrinkled} and~\cite{Williams} ) to a Morse $2$-function, which are also the changes
possible in a Cerf graphic.  These are:

\begin{enumerate}

\item (swallowtail) Given a $k$-fold, a $(k\pm 1$)-fold can be
  introduced as in Figure~\ref{F:Moves}a.  The swallowtail can be removed if it is
  isolated as in Figure~\ref{F:Moves}a.

\item (eye)  A canceling pair of $k$- and $(k+1)$-folds can be
  introduce, or removed, as in Figure~\ref{F:Moves}b.

\item (merge) If two $(k, k+1)$-cusps appear as in Figure~\ref{F:Moves}c, then they
  can be {\it merged} to form two separate folds.  The reverse (an
  {\it unmerge} or {\it cancellation}) can occur if the $(k+1)$-fold
  cancels the $k$-fold, in the sense that at some given time, the
  $(k+1)$-handle cancels the $k$-handle.

\end{enumerate}

\begin{theorem}[Uniqueness]  Given two fiber connected, indefinite, Morse
$2$-functions $X^n \to \Sigma$ which agree on $\partial X$, then the
generic homotopy between them, using  regular homotopies  and the three local
changes above, can also be chosen to be indefinite at all times if
$n>3$, and to have connected fibers throughout.
\end{theorem}

The indefinite part of this theorem for $n=4$ was proved independently
and earlier by Jonathan Williams \cite{Williams}.  He also showed that
for the $S^2$ case, it is possible to choose a map with only one
connected fold, with a genus $g$ surface as fiber on one side and a
genus $g+1$ fiber on the other side.

In the proof, removing definite folds, which may be introduced as
swallowtails or eyes, involves the use of codimension two
singularities, in particular the {\it butterfly} and the {\it elliptic
  umbilic}. To keep fibers connected involves old fashioned bordism
arguments  and careful control of regular homotopies .

Much of the original motivation for this work, in the $4$--dimensional setting, comes from the work of Auroux, Donaldson and Katzarkov~\cite{ADK} on singular Lefschetz fibrations, which became known as broken (Lefschetz) fibrations, the following work of Perutz~\cite{PerutzI, PerutzII} setting up a Floer theoretic framework in which to define invariants in terms of broken fibrations and near-symplectic forms, and the work of Lekili~\cite{LekiliWrinkled}, Baykur~\cite{BaykurExist} and Williams~\cite{Williams} relating homotopic broken Lefschetz fibrations to indefinite Morse $2$--functions. Maps to $S^2$ mixing Lefschetz singularities with indefinite Morse $2$--function singularities were called {\em wrinkled fibrations} in this earlier work, and were called {\em purely wrinkled fibrations} in the absence of any Lefschetz singularities. 

Our contribution of fiber-connectedness is important if one wants to work with near-symplectic fibrations and keep the regular fibers symplectic through homotopies. We also stress that our existence and uniqueness results work just as well over surfaces with boundary, another new feature, and that this is potentially critical to developing invariants for cobordisms between circle-valued Morse functions.

\textbf{Proof of existence, Theorem~\ref{T:Existence}, when $\Sigma=S^2$:} This proof is quite different from that in
\cite{Saeki} but suggestive of some of the work in \cite{GayKirbyMorse2Fs}.  Choose any
codimension $2$, orientable, closed, submanifold $Z^{n-2}$ of $X$ with
the property that its normal bundle $\nu$ is trivial, and that the
projection of the normal circle bundle $Z \times S^1 \to S^1$ extends
to a map $g:X-Z \to S^1$ for some trivialization of $\nu$.  Note that
such a $Z$ always exists.  For consider the case of $S^n = (S^{n-2}
\times B^2) \cup (S^1 \times B^{n-1})$ where $S^{n-2} \times 0 = Z$.
Then write $X = X \# S^n = (X^n -B^n) \cup (S^n -B^n)$ and extend the
map on $S^n - B^n$ over $X^n -B^n$ by essentially mapping $X^n - B^n$
to the neighborhood of a point in $S^1$.

A special case, which is useful to start with, is that of an open book
with connected fiber. (Ka Choi in~\cite{Choi} gives some explicit examples of fiber connected, indefinite Morse $2$-functions $S^4 \to S^2$ in which the fiber over the north pole is a $2$-sphere which is a fibered $2$-knot in $S^4$).  $X^n$ is described as $Z \times B^2) \cup
(F^{n-1} \tilde{\times} S^1)$ where $F\tilde{\times} S^1 \to S^1$ is a
bundle with fiber $F$ and monodromy $\mu: F \to F$ fixing $\partial
F$.  The boundary of $F$ is $Z$ and $(z,1,\theta)$ is identified with
$(z,\theta)$. 

We begin the definition of $p:X \to S^2$ by sending $Z \times B^2$ to
the northern hemisphere by projection on the second coordinate;  $Z
\times 0$ goes to the north pole.  Now parameterize longitudes in the
southern hemisphere by $[-1,n-1]$ with $-1$ at the south pole and
$n-1$ at the equator.  Next pick a Morse function $f:F^{n-1} \to
[-1,n-1]$ with the critical points of index $k$ being mapped to
distinct values in a small neighborhood of $k \in [-1,n-1]$, and with
only one critical point of index $0$.

If the monodromy $\mu$ is the identity, we can then extend $p$ over $Z
\times S^1$ by $p(y,\theta) = (f(z), \theta)$. where the latter pair
lies on the longitude at $\theta$.  Then the folds will be latitudes
of the southern hemisphere at various heights depending on the
critical values of $f$.

This is not yet an indefinite Morse $2$-function because of the
minimum of $f$ which generates a $0$-fold at, say, the antarctic
circle.  This can be removed, but first we deal with the case of
non-trivial monodromy $\mu$ and a generalization to arbitrary $X$.

\begin{figure} 
\begin{center}
\includegraphics{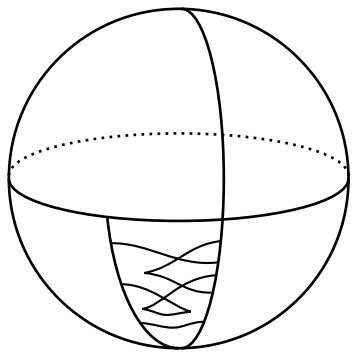}
\caption{Constructing a $S^2$--valued Morse $2$--function from an open book decomposition.} \label{F:OpenBook}
\end{center}
\end{figure}

For non-trivial $\mu$, consider the Morse function on $Z$ given by the
composition $f \circ \mu : Z \to [-1,n-1]$.  There is a nice homotopy
between $f$ and $f\mu$ producing a {\it Cerf graphic}, and this graphic
can be spliced into the southern longitudes near Greenwich zero, as
suggested in Figure~\ref{F:OpenBook}.

\begin{figure} 
\begin{center}
\includegraphics{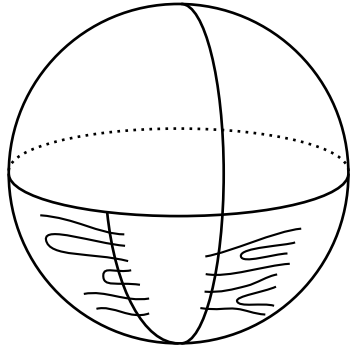}
\caption{The more general case where the ``open book'' has isolated singular pages.} \label{F:BrokenOpen}
\end{center}
\end{figure}

For arbitrary $X$ with $Z^{n-2}$ and $g:X-Z \to S^1$ with no local
minima or maxima, we begin with the fiber $F = g^{-1}(0)$ for $0$ a
regular value.  Assume $f:F \to [-1,n-1]$ is a Morse function as
before, and use it to map $F$ to the $0$-longitude.  This extends
trivially to a map $p:F \times [-\epsilon, \epsilon] \to S^2$.  When a
critical point of $g$ of index $k$ is encountered, the folds will
change by the addition of a $k$-handle which can be drawn as a new
$k$-fold with vertical tangent, as in Figure~\ref{F:BrokenOpen}.

In this way, we progress past all of the critical points of $g$ and
are left  with a monodromy $\mu$ which is used to provide a Cerf
graphic filling in the remaining gap, as in the case of an open book
above. 

\begin{figure} 
\labellist
\tiny\hair 2pt
\pinlabel $0$ [b] at 22 40
\pinlabel $1$ [b] at 65 39
\pinlabel $0$ [r] at 71 22
\pinlabel $1$ [b] at 108 36
\pinlabel $0$ [r] at 97 22
\pinlabel $0$ [r] at 135 22
\pinlabel $1$ [b] at 151 16
\pinlabel $1$ [b] at 194 11

\small
\pinlabel (a) at 22 1
\pinlabel (b) at 64 1
\pinlabel (c) at 108 1
\pinlabel (d) at 151 1
\pinlabel (e) at 194 1

\endlabellist

\begin{center}
\includegraphics{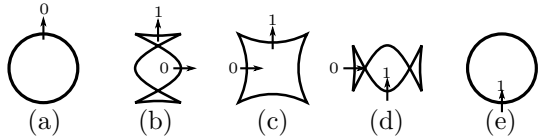}
\caption{Flipping a $0$--fold to a $1$--fold.} \label{F:ZeroToOne}
\end{center}
\end{figure}

What remains for the proof of existence is to remove the $0$-fold at
the antartic circle.  In a coordinate patch around the south pole, we
see a $0$-fold looking out, or an $(n-1)$-fold looking in.  First we
introduce a pair of swallowtails, as drawn in Figure~\ref{F:ZeroToOne}b.  The
$0$-folds can then be pulled backwards, simply because the birth of a
minimum can always occur earlier.  In figure~\ref{F:ZeroToOne}c, we can pull both
$1$-folds  across each other because in a Cerf graphic one can pull a
$2$-fold below a $1$-fold.  The last step, from Figure~\ref{F:ZeroToOne}d to Figure~\ref{F:ZeroToOne}e
is to cancel the two swallowtails leaving a $1$-fold pointing in.

The fiber at the south pole is now $S^1 \times S^{n-3}$.  In a
non-obvious way, this fiber is glued in by the map which rotates the
circle factor as one traverses the antartic circle.  This phenomenon
was first discovered in \cite{ADK}, and is covered in detail in
\cite{GayKirbyMorse2Fs}.
$\Box$

Note that, at the end of the construction above, we have a fiber of genus $1$ over the south pole, surrounded by a circle of $2$--folds pointing out, with $2$--sphere fibers outside that. This is the example where a nontrivial element of $\pi_1(\mathit{Diffeo}_+(S^1 \times S^1))$ appears. Note that, whenever a genus $1$ fiber $F$ occurs, we can remove a neighborhood of $F$ and glue it back in using a nontrivial element of $\pi_1(\mathit{Diffeo}_+(S^1 \times S^1))$ and, usually, change the $4$--manifold.  This is the usual {\em logarithmic transform} operation . However, this operation does not change the local data (the attaching $S^1$'s and $S^0$'s in one fiber from each region). 

\begin{figure} 
\begin{center}
\includegraphics{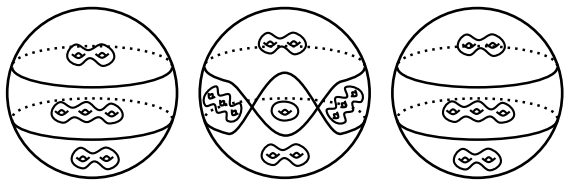}
\caption{Monodromy in non-simply connected regions and $\pi_1(\mathit{Diffeo}_+(S^1 \times S^1))$.} \label{F:Monodromy}
\end{center}
\end{figure}

Here is an interesting example where the indeterminacy associated with $\pi_1(\mathit{Diffeo}_+(S^1 \times S^1))$ is traded for the indeterminacy associated with a non-simply connected region in $S^2 - F$ (where $F$ is the image of the fold locus in $S^2$): The example is illustrated in Figure~\ref{F:Monodromy}. Start with $p : X^4 \to S^2$ in which $F$ is two circles, the tropics of Cancer and Capricorn. Over the equator the fiber $F$ has genus $3$ and over the poles the fiber has genus $2$, so that the folds are $2$--folds going towards the poles. Choose $p$ so that the attaching circles for the $2$--folds in $F$ are disjoint. This means that, over a short section of the equator, we can pull one $2$--fold below the other, producing a bigon over with the fiber has genus $1$. Now remove the inverse image of a disk in this region and glue back in by a nontrivial element of $\pi_1(\mathit{Diffeo}_+(S^1 \times S^1))$, and then pull the $2$--folds apart again. The resulting picture looks just like the beginning picture, but the $4$--manifold  may have  changed because the monodromy has changed around the annular equatorial region.

\end{document}